\documentclass[12pt]{article}
\newcommand{\qed}{\ifmmode$\Box$\else{\unskip\nobreak\hfil
\penalty50\hskip1em\null\nobreak\hfil$\Box$
\parfillskip=0pt\finalhyphendemerits=0\endgraf}\fi}
\usepackage{amsmath}
\usepackage{latexsym}
\usepackage{amsfonts}
\begin{document}
\renewcommand{\theequation}{\thesection.\arabic{equation}}
\newtheorem{Pa}{Paper}[section]
\newtheorem{Tm}[Pa]{{\bf Theorem}}
\newtheorem{La}[Pa]{{\bf Lemma}}
\newtheorem{Cy}[Pa]{{\bf Corollary}}
\newtheorem{Rk}[Pa]{{\bf Remark}}
\newtheorem{Pn}[Pa]{{\bf Proposition}}
\newtheorem{Ee}[Pa]{{\bf Example}}
\newtheorem{Dn}[Pa]{{\bf Definition}}
\newtheorem{I}[Pa]{{\bf}}
\oddsidemargin 0in
\topmargin -0.5in
\textwidth 16.5truecm
\textheight 23truecm
\parindent =0.cm
\textwidth 16.5truecm
\textheight 23truecm
\begin{center}
{\bf \Large On the structure of Banach algebras associated with automorphisms. 2.
Operators with measurable coefficients}\\
\bigskip
{\large A. Lebedev \\ }
\begin{center} {Institute of Theoretical Physics\\
University in Bia{\l}ystok
\\Lipowa 41, 15-424 Bia{\l}ystok, Poland}
\end{center}
\end{center}
\bigskip
\begin{abstract}
In the present paper we continue the
study of the structure of a Banach algebra $B(A, T_g )$ generated
by a certain Banach algebra $A$ of operators acting in a Banach space $D$ and a group
$\{ T_g \}_{g \in G}$ of isometries of $D$  such that $T_g A T^{-1}_g = A$.
We investigate the interrelations
between the existence of the  expectation of $B(A, T_g )$ onto $A$, metrical freedom of
the automorphisms of $A$ induced by $T_g$ and  the dual action of the group $G$ on $B(A, T_g )$.
The  results obtained are applied to the description of the structure of Banach algebras generated
by 'weighted composition operators' acting in Lebesgue spaces.
\end{abstract}



\tableofcontents
\section{Introduction}
\setcounter{equation}{0} This article should be considered as a
'measurable counterpart' of \cite{Leb}. As in \cite {Leb} the
principal  object under consideration here   is a Banach algebra
$B(A, T_g )$ generated by a certain Banach algebra $A$ of
operators acting in a Banach space $D$ and a group $\{ T_g \}_{g
\in G}$ of isometries of $D$ (a representation $g \to T_g$ of a
discrete group $G$) such that
\begin{equation}
\label{ae1.1}
T_g A T^{-1}_g = A, \hspace{10mm} g\in G
\end{equation}
which means that $T_g$ generates the automorphism ${\hat T}_g$ of $A$ given by
\begin{equation}
\label{ae1.2}
{\hat T}_g (a) = T_g a T^{-1}_g , \hspace{10mm} a\in A.
\end{equation}
In \cite{Leb} we obtained some principle characteristics
of the structure of such  algebras  and also considered a
number of examples where the role of $A$ was played
by algebras of continuous operator valued functions. In the present article we continue
this investigation and present the results on the structure of the corresponding
algebras in the situation when as $A$ are  taken algebras of measurable operator
valued functions acting in Lebesgue spaces.
Therefore roughly speaking the material of the paper
describes a passage from the 'topological'
picture given in \cite{Leb} to a 'measurable' one.\\

 We recall that in the Hilbert space situation (that is in the    $C^\ast -$algebra theory)
 the analogous objects
 are closely related to the crossed products
(see, for example \cite{Ped}) and  description of their structure
is the theme of numerous investigations. In particular, Landstad
\cite{Land} presented the necessary and sufficient conditions (in
terms of {\em duality theory}) for a $C^\ast -$algebra to be
isomorphic to a crossed product (of an
 algebra and a locally compact group of automorphisms).
In the case of a discrete group in \cite{AnLeb},
Chapter 2 there were found the conditions for a $C^\ast -$algebra to be isomorphic to
a crossed product in terms of the group
action (the so called {\em topologically free action} (see
\ref{1.6})) and also in
terms of satisfaction of a certain  inequality ({\em property (*)}
(\ref{e1.3}))
guaranteeing the existence of the  expectation  of the algebra $B(A, T_g )$ onto the
algebra $A$ (see (\ref{e1.4}), (\ref{e1.5})).\\

In \cite{Leb} we have shown that
 the mentioned properties
(topologically free action, property (*) and dual action of the
group) play an exceptional role in the general Banach space
situation as well. By means of these properties there were
established a number of results describing the structure of
 of $B(A, T_g )$ up to isomorphism.\\

In this paper we show (in Section 2) that in the 'measurable' situation the natural
substitute for the topologically free action is the {\em metrically free} action (see
\ref{a1.9}). We investigate the interrelation between these notions and in
particular find out that from a ceratin point of view they are equivalent.
This enables us to transfer the main structural results of \cite{Leb}
from the 'topological' to the 'measurable' situation.\\

Since in the general Banach space situation we do not have a
universal object like the crossed product in the Hilbert space situation to
describe the structure of
$B(A, T_g )$ we have to specify the algebra $A$ and the isometries $T_g$.
Sections 3-5 are devoted
to the applications of the results obtained in Section 2 and in \cite{Leb}
to the description of the structure of
concrete Banach algebras associated with automorphisms,
namely the algebras generated by
'weighted composition operators' acting   in Lebesgue  spaces.\\

We establish a number of isomorphism results for the algebras investigated and in addition
find out that the arising algebras are in a way qualitatively different. In particular when
considering the operators in $ L^\infty _\mu (\Omega , E)$ and in
$ L^1 _\mu (\Omega , E)$ we {\em can}
calculate their norms  (see Theorem \ref{a2.21} and \ref{a2.30}) while for the operators in
$L^p _\mu (\Omega , E), \hspace{2mm} 1<p < \infty$ we have nothing like this.
Moreover  to  obtain the  isomorphism
theorems  for the algebras
$B(A, T_g )$ in   $ L^\infty _\mu (\Omega , E)$, $ L^1 _\mu (\Omega , E)$ we
do not need any information on the structure of the group of operators generating automorphisms
while this structure (namely the amenability of the group $G$) is vital when we are
investigating the operators
in $L^p _\mu (\Omega , E), \hspace{2mm} 1<p < \infty $
(see Theorem \ref{a2.14} and Remark \ref{a2.33}).\\

To make the presentation selfcontained we have
 to recall a number of notions and results from \cite{Leb}.
\begin{I}
\label{1.3}
{\bf Property (*).}
\em
It was shown in \cite {Leb} that one of
 the most important properties of the algebra  $B(A, T_g )$ in the
presence of which one can obtain the deep and fruitful theory of the subject is the
next {\em property (*)}:\\

for any finite sum $b = \sum a_g T_g , \hspace{3mm}a_g \in A$ the following inequality holds
\begin{equation}
\label{e1.3}
\Vert b \Vert = \Vert  \sum a_g T_g \Vert \ge \Vert a_e \Vert ,
\end{equation}
where $e$ is the identity of the group $G$.
\end{I}
If an algebra $B(A, T_g )$ possesses the property (*)
then for every $g_0 \in G$ there is correctly defined
the mapping
\begin{equation}
\label{e1.4}
N_{g_0} : \sum a_g T_g \to a_{g_0}
\end{equation}
which can be extended up to the mapping
\begin{equation}
\label{e1.5}
N_{g_0} : B(A, T_g ) \to A.
\end{equation}
 \begin{I}
\label{1.4}
{\bf Property (**).}
\em
One more important property is the following.\\
We shall say that an algebra $B(A, T_g )$ which possesses the property (*)
also possesses  the {\em property (**)}
if
\begin{equation}
\label{e1.6}
B(A, T_g ) \ni b = 0 \hspace{2mm}{\rm iff}  \hspace{2mm}N_g (b) = 0 \hspace{2mm} {\mbox {\rm for
 every}}  \hspace{2mm}g \in G
\end{equation}
where $N_g$ is the mapping introduced above.
\end{I}
In fact the presence of the properties (*) and (**) makes it possible to 'reestablish' an
element $b \in B(A, T_g )$ via its 'Fourier' coefficients $N_g (b), \hspace{2mm}g\in G$
and we shall find out further that in many reasonable situations this 'reestablishing'
can be carried out successfully.\\

If $A$ is a $C^\ast -$algebra of operators containing the identity and
acting in a Hilbert space $H$ and $\{ T_g \}_{g \in G}$
is a unitary representation of a group $G$ in $H$ then the $C^\ast -$algebra generated
by  $A$ and  $\{ T_g \}_{g \in G}$ will be denoted by $C^\ast (A, T_g )$.

In the $C^\ast -$algebra situation we have (see \cite{AnLeb}, Theorems 12.8 and 12.4):\\

{\em if $G$ is a discrete amenable group and $C^\ast (A, T_g )$
possesses the property (*)
then $C^\ast (A, T_g )$ possesses the property (**) as well. }\\

The main reason why in the $C^\ast -$algebra case the
property (**) (\ref{e1.6}) follows from the property (*) is that\\

 {\em the presence of the property (*) implies}  (\cite{AnLeb}, Theorem 12.8)
$$
C^\ast (A, T_g ) \cong A \times_{\hat T} G
$$
where by $ A \times_{\hat T} G$ we denote the cross-product of the algebra $A$ by the group
$\{ {\hat T}_g \}_{g \in G}$ of its automorphisms (here $G$ is considered as a discrete group).\\
Since in a Banach space case we do not have anything like the
isomorphism mentioned above we have to check the property (**)
even when $B(A, T_g )$ possesses the property (*). In a general
situation (that is for an {\em arbitrary} discrete group of
isometries $\{ {\hat T}_g \}_{g \in G}$ with $T_g A T^{-1}_g = A$)
the verification of the property (**) may be very sophisticated.
 The next Theorem \ref{1.13} (proved in \cite{Leb})
shows that in the case of a {\em locally compact commutative}
group $G$ and under a special assumption (which as it will be seen
later is in fact rather common) the algebra $B(A, T_g )$ possesses
the properties (*) and (**) simultaneously.
\begin{Tm}
\label{1.13}
Let $G$ be a locally compact commutative group. If for any finite set $F\subset G$ and
any character $\chi \in {\hat G}$ there is satisfied the equality
\begin{equation}
\label{e1.36}
\Vert \sum_{g\in F} a_g T_g \Vert = \Vert \sum_{g\in F} a_g  \chi (g) T_g \Vert
\end{equation}
then the algebra $B(A, T_g )$ possesses the properties (*) and (**).
\end{Tm}
In \cite{Leb} we have also established  a close  relation between the property (*)
and the so-called {\em topological freedom} of the action of the group of automorphisms
$\{ {\hat T}_g \}$.  So let as recall the latter notion.
\begin{I}
\label{1.6}
{\bf Topologically free action.}
\em
Observe first that if
 $\{ T_g \}_{g\in G}$ is a group of isometries satisfying (\ref{ae1.1}) then evidently
\begin{equation}
\label{center}
T_g {\mathbb Z}(A) T_g^{-1} = {\mathbb Z}(A)
\end{equation}
 where ${\mathbb Z}(A)$ is the center of $A$.\\
 Let $A$ be a certain Banach algebra isomorphic to $C(X, B)$ where $X$ is a
completely regular space and $B$ is a Banach algebra then
\begin{equation}
\label{zenter1}
{\mathbb Z}(A) = C(X, {\mathbb Z}(B)).
\end{equation}
Henceforth in this subsection  we  confine ourselves to the case
\begin{equation}
\label{center5}
{\mathbb Z}(B) =  \{ cI \}
\end{equation}
The reason justifying this choice was discussed in \cite{Leb}.
 Obviously  if $B = L(E)$ is the Banach algebra of all
linear bounded operators acting in a Banach space $E$
 then  (\ref{center5}) is satisfied.\\
So let $A\subset L(D)$ be a Banach algebra of operators isomorphic
to $C(X, L(E))$ where $X$ is a certain completely regular space
and $E$ and $D$ are Banach spaces (thus ${\mathbb Z}(A) \cong
C(X)$). Let $\{ T_g \}_{g\in G}$ be a group of isometries
satisfying (\ref{ae1.1}).  According to (\ref{center}) the
automorphisms ${\hat T}_g$  (\ref{ae1.2}) preserve the center and
henceforth we assume that their action on the center is given by
\begin{equation}
\label{e1.7}
[{\hat T}_g (z)](x) = z(t_g^{-1} (x)), \hspace{5mm}z\in {\mathbb Z}(A), \hspace{2mm} x\in   X.
\end{equation}
where $t_g : X \to X$ are some homeomorphisms of $ X $.\\

Denote by $X_g , \hspace{2mm} g\in G  $ the set
\begin{equation}
\label{e1.8}
X_g = \{ x \in X : t_g (x) = x \}.
\end{equation}
We say that the group $G$ acts {\em topologically freely} on $A$
by automorphisms ${\hat T}_g$ (or on $X$ by homeomorphisms $t_g$
mentioned in (\ref{e1.7}) if for any $g\in G,\hspace{2mm} g\ne e$
 the set
$ X_{g} $
has an empty interior.\\

One can observe that  $G$ acts  topologically freely iff   for any
finite set $\{ g_1 , ... , g_n \} \subset G \hspace{2mm}
(g_i \ne e)$ the set $[ \cup^n _{i = 1} X_{g_i} ]$ has an empty interior.\\

Just as in \cite{AnLeb}, 12.13 and 12.13' it can be shown that the
foregoing definition is equivalent to the next one: the action of
$G$ is said to be topologically free if for any finite set  $\{
g_1 , ... , g_k \} \subset G$ and a non empty open set $U\subset
X$ there exists a point $x\in U$ such that all the points $t_{g_i}
(x), \hspace{2mm}
i=1, ... , k$ are distinct.\\

Since $X$ is Hausdorff the latter definition is also equivalent to the following.\\
The action of $G$ is said to be topologically free if for any
finite set $\{ g_1 , ... , g_k \} \subset G$ and a non empty open
set $U\subset X$ there exists a non empty open set $V\subset U$
such that
\begin{equation}
\label{e1.6'}
t_{g_i} (V) \cap t_{g_j} (V) =\emptyset \hspace{5mm}i,j \in \overline{1,k},\hspace{2mm}i\ne j.
\end{equation}
\end{I}
In \cite{Leb} we proved the following
\begin{Tm}
\label{1.7}
If\ \ $G$ acts topologically freely then $B(A, T_g )$ possesses the property (*).
\end{Tm}
\begin{I}
\label{1.17}
\em
{ \bf Regular representation of an algebra and a group of automorphisms. }
Let us also recall one more algebra (examined in \cite{Leb})
 where the properties (*) and (**) can be checked easily
--- {\em the regular representation of an algebra $A$ and a group of automorphisms
$\{ {\hat T}_g \}_{g\in G}$. }\\
Namely let $A\subset L(D)$ be a certain Banach algebra and $\{ {\hat T}_g \}_{g\in G}$
be a certain
group of its automorphisms ($G$ is an {\em arbitrary} group that is {\em not necessarily
commutative}).\\
Denote by $H$ any of the spaces $l^p (G, D), 1\le p \le \infty$ or $l_0 (G, D)$
(here $l_0 (G, D)$ is the space of vector valued functions on $G$ having values in $D$
and tending to zero at infinity (with the sup-norm)).\\
Set the operators $V_{g_0} : H \to H$ by the formula
\begin{equation}
\label{e1.51}
(V_{g_0} \xi )(g) = \xi (gg_0 ), \hspace{5mm} g, g_0 \in G
\end{equation}
and consider the algebra ${\bar A} \subset L(H)$ isomorphic (as a Banach algebra)
 to $A$ and given by
\begin{equation}
\label{e1.52}
({\bar a}\xi )(g) = {\hat T}_g (a) \xi (g), \hspace{5mm} a \in A.
\end{equation}
Routine computation shows that with this notation we have
$$
V_g {\bar a} V^{-1}_g = {\overline{{\hat T}_g (a)}}
$$
which in view of the isomorphism between $A$ and $\bar A$ means that the operators $V_g , g\in G$
given by (\ref{e1.51}) generate the automorphisms ${\hat T}_g$ of $\bar A$.\\

The algebra $B({\bar A}, V_g )\subset L(H)$ is called the {\em (right) regular representation}
corresponding to the algebra $A$ and the group of automorphisms $\{ {\hat T}_g \}_{g\in G}$
(in fact we have the series of representations depending on the type of the space
$H$ chosen).
\end{I}
In \cite{Leb} we have proved  that \\

{\em the algebra $B({\bar A}, V_g )$ possesses the properties (*),  (**) and  (\ref{e1.36})
(for every $H$ considered)}.\\

\section{Metrically free action and topologically free action}
\setcounter{equation}{0}

\begin{I}
\label{a1.8}
\em
Let $(\Omega , \mu )$ be a space with a $\sigma -$additive $\sigma -$finite measure $\mu $,
$H$ be a certain Banach space and $L^p_\mu (\Omega , H ), \ \ 1\le p \le \infty$ be the spaces of
(equivalence classes) of measurable functions $f : \Omega \to H$ bounded with respect to the norms
$$
\Vert f \Vert = \left[ \int_\Omega \vert f(x)\vert^p \, d\mu (x) \right] ^{\frac{1}{p}}, \ \ 1\le p< \infty ,
$$
$$
\Vert f \Vert = {\rm esssup}_{\Omega} \vert f \vert ,\ \ \ \ \ \ \ \ \ \ \ \ \ \ p=\infty
$$
where $\vert \cdot \vert$ is the norm in $H$ (for details see
for example Dunford, Schwarts \cite{DunSchw}).\\
Consider an algebra $A \subset L(D)$ isomorphic to $L^\infty _\mu (\Omega , L(E))$
where $D$ and $E$ are
Banach spaces. If $\{ T_g \}_{g\in G}$ is a group of isometries of $D$
satisfying (\ref{ae1.1}) then
the automorphisms ${\hat T}_g$ (given by (\ref{ae1.2})) generate the mappings
$\alpha_g : \Sigma \to \Sigma $
($\Sigma $ is the set of (equivalence classes) of measurable subsets of $\Omega$)
defined in the following way.\\

Observe that for the algebra considered the center
${\mathbb Z}(A) \cong L^\infty _\mu (\Omega) $. Let $\chi _\Delta $ be the element of
${\mathbb Z}(A)$ corresponding to the characteristic function  $\chi _\Delta (\omega )$
of a certain set $\Delta \in \Sigma$.\\
Since ${\chi _\Delta}^2 =\chi _\Delta  $ it follows that ${\hat T}_g (\chi _\Delta   )$
is a projection belonging to ${\mathbb Z}(A)$ (non zero iff $\chi _\Delta  \neq 0$)
and we have
\begin{equation}
\label{ae1.15}
{\hat T}_g (\chi _\Delta   ) = \chi_{\tilde \Delta}
\end{equation}
for some ${\tilde \Delta} \in \Sigma$.\\
We set
\begin{equation}
\label{ae1.16}
\alpha_g (\Delta ) ={\tilde \Delta} .
\end{equation}

The substitution for the topologically free action of $G$  (see \ref{1.6})
in the
 situation under consideration is the
so-called metrically free action. Here it is.
\end{I}

\begin{I}
\label{a1.9}
\em
We say that the group $G$ acts {\em metrically freely} on $A$ (considered in \ref{a1.8})
by automorphisms ${\hat T}_g$ (or on $\Sigma $ by the mappings $\alpha _g$) if
for any finite set $\{ g_1 , . . . , g_k \} \subset G$ and any $\Delta \in \Sigma$ with
$\mu (\Delta )> 0$ there exists a set $\Delta ^\prime \in \Sigma$ such that

(i) $\mu (\Delta^\prime ) > 0$,

(ii) $\Delta^\prime \subset \Delta$,

(iii) $\mu (\alpha_{g_i} (\Delta^\prime ) \cap \alpha_{g_j} (\Delta^\prime )) = 0, \ \ \
i,j\in\overline{1,k}, \ i\neq j$.
\end{I}
\begin{Rk}
\label{a1.10} \em It is worth mentioning that from a certain point
of view the notion of the metrically free action of $G$ just
introduced 'coincides' with the notion of the topologically free
action.
\end{Rk}

Indeed.\\
The algebra $L^\infty _\mu (\Omega )$ is a commutative $C^\ast -$algebra
(with the natural involution). Let $M$ be its maximal ideal space then
$$
L^\infty _\mu (\Omega ) {\check{\cong}} C(M)
$$
where the isomorphism is established by means of the Gelfand transform.
Under this isomorphism the characteristic function
$\chi _\Delta (\omega ) \in L^\infty _\mu (\Omega ) $
goes into a certain function ${\check\chi }_\Delta \in C(M)$ which
(since ${\chi _\Delta}^2 = \chi _\Delta$) is also a characteristic function of a certain set
${\check \Delta} \subset M$.
If $\mu (\Delta )>0$ (that is $\chi _\Delta \neq 0$) then
${\check\chi} _\Delta = \chi_{\check \Delta}  \neq 0$ which means that
${\check \Delta}$ is a non empty open set.\\
Now observe that\\

{\em
 for any open non empty set $U\subset M$ there exists a set $\Delta \in \Sigma$
with the properties
\begin{equation}
\label{ae1.17}
\mu (\Delta ) > 0 \ \ and \ \ {\check \Delta}\subset U.
\end{equation}
}\\

{\bf Proof:} Note first that the measure $\mu $ generates a certain measure
${\check \mu}$ on $M$ with the properties
\begin{equation}
\label{ae1.18}
\int_\Omega a \, d\mu = \int_M {\check a} \, d{\check \mu}
\end{equation}
and
\begin{equation}
\label{ae1.19}
{\rm supp}\, {\check \mu} = M.
\end{equation}

Indeed let
$$
\Omega = \cup^\infty _{i = 1} \Omega _i
$$
where $\Omega_i \in \Sigma , \ \  \Omega_i \cap \Omega_j = \emptyset , \ \ i\neq j$
and $\mu (\Omega_i ) < \infty$.\\
Then the mapping
$$
{\check a} \to \int_{\Omega_i } a\, d\mu
$$
is a linear positive functional. Thus there exists a regular Borel measure
${\check \mu}_i$ on $M$
such that
$$
\int_M {\check a}\, d\mu_i = \int_{\Omega_i } a\, d\mu .
$$
In particular
$$
\mu (\Omega_i ) = \int_\Omega \chi_{\Omega_i}\, d\mu =
\mu_i (M) = \mu_i ({\check \Omega}_i )
$$
and the condition $\Omega_i \cap \Omega_j =\emptyset, \ i\neq j$
implies
$$
\mu_i ({\check \Omega}_j )= 0 \ \  {\rm and} \ \ \mu_j ({\check \Omega}_i )= 0.
$$
Now set
$$
{\check \mu} = \sum_i \mu_i .
$$
By the construction ${\check \mu} $ satisfies (\ref{ae1.18}), and since
$\cup^\infty _{i=1} \Omega_i =\Omega$ (\ref{ae1.19}) holds as well.\\
If $U=M$ then  (as ${\check \Omega} = M$) we can take in (\ref{ae1.17})
$\Delta = \Omega $.\\
Now let $U$ be any open non void subset of $M \ \ (U\neq M)$.
Since $M$ is normal there exists an open non void set
$U_1$ with the properties
$$
U_1 \subset \overline{U_1} \subset U \ \ {\rm and} \ \ \overline{U_1} \cap
[M\setminus U] =\emptyset .
$$
By the Urysohn Lemma there exist functions
${\check \varphi}_1 , {\check \varphi}_2 \in C(M)$
with the properties

(1) $0\le {\check \varphi}_i (x) \le 1\ \  i=1,2$,

(2) ${\check \varphi}_1 (x^\prime )= 1$ for some $x^\prime \in U_1$,

(3) ${\check \varphi}_1 (M\setminus U_1 ) = 0 $,

(4) ${\check \varphi}_2 ( \overline{U_1}) = 0$,

(5) ${\check \varphi}_2 (M\setminus U ) = 1 $.\\
By (2) and (5) and (\ref{ae1.19}) it follows that
\begin{equation}
\label{ae1.20}
\int_M {\check \varphi}_i (x) \, d{\check \mu} > 0, \ \ \ i=1,2.
\end{equation}
By (3) and (4) we have
\begin{equation}
\label{ae1.21}
{\check \varphi}_1 \cdot {\check \varphi}_2 = 0.
\end{equation}
In turn (\ref{ae1.20}), (\ref{ae1.18}) and (1) imply
\begin{equation}
\label{ae1.22}
\int_\Omega \varphi_i \, d\mu = \int_{\Delta_i}\varphi_i \, d\mu > 0 \ \ i=1,2
\end{equation}
where
\begin{equation}
\label{ae1.23}
\Delta_i = \{ \omega \in \Omega : \varphi_i (\omega ) > 0 \} .
\end{equation}
Note that in particular (\ref{ae1.22}) implies
\begin{equation}
\label{ae1.24}
\mu (\Delta_i ) >0 \ \ i=1,2.
\end{equation}
In view of (\ref{ae1.21}) we have
\begin{equation}
\label{ae1.25}
\varphi_1 \cdot \varphi_2 = 0
\end{equation}
which means that
$$
\mu (\Delta_1 \cap \Delta_2 ) =0
$$
thus
$$
\chi_{\Delta_1}\cdot \chi_{\Delta_2} =0
$$
and it follows that
\begin{equation}
\label{ae1.26}
{\check \chi}_{\Delta_1} \cdot {\check \chi}_{\Delta_2} =
\chi_{{\check \Delta}_1} \cdot  \chi_{{\check \Delta}_2} = 0.
\end{equation}
Observe also that by the definition of $\Delta_i$ we have
\begin{equation}
\label{ae1.27}
\chi_{\Delta_i}\cdot \varphi_i =\varphi_i \ \ i=1,2
\end{equation}
and consequently
\begin{equation}
\label{ae1.28}
{\check \chi}_{\Delta_i}\cdot {\check\varphi}_i = {\check\varphi}_i \ \ i=1,2
\end{equation}
which in turn along with (\ref{ae1.26}) implies
\begin{equation}
\label{ae1.29}
\chi_{{\check \Delta}_1}\cdot {\check\varphi}_2 =
{\check \chi}_{\Delta_1}\cdot {\check \chi}_{\Delta_2}\cdot {\check\varphi}_2
\end{equation}
and (\ref{ae1.29}) means in particular that
$$
{\check \Delta}_1 \cap (M\setminus U) = \emptyset
$$
that is
\begin{equation}
\label{ae1.30}
{\check \Delta}_1 \subset U.
\end{equation}
And in view of (\ref{ae1.24}) and (\ref{ae1.30}) we conclude that
$\Delta_1$ is the desired set.
\mbox{}\qed\mbox{}\\

({\bf Remark.} A slight improvement of the foregoing argument can
show that there exists a set $\Delta \in \Sigma $ with the
properties
$$
\mu (\Delta ) >0 \ \ \  {\rm and} \ \ \ {\check \Delta} = U
$$
thus establishing the correspondence between the elements of $\Sigma $
having non zero measure and non empty open subsets of $M$.)\\

Now we are ready to establish the 'coincidence' of the
topologically
and metrically free actions. Namely we shall prove that\\

{\em
the metrically free action of the automorphisms ${\hat T}_g$ on
$L^\infty _\mu (\Omega )$ corresponds to the topologically free
action of the automorphisms
${\check T}_g$ on $C(M)$ (induced by the automorphisms
${\hat T}_g$ and the isomorphism
$L^\infty _\mu (\Omega ) {\check {\cong}} C(M)$).
}

{\bf Proof:}

{\bf I.} Let $\{ {\hat T}_g \}_{g\in G}$ act metrically freely. \\
Consider any finite set $\{ g_1 , . . . , g_k \} \subset G$, any open set
$U\subset M$ and the set $\Delta \in \Sigma$ defined by (\ref{ae1.17}).\\
Now let  $\Delta ^\prime \in \Sigma$ be the set mentioned in the definition of
the metrically free action (in \ref{a1.9}). The condition (ii) in \ref{a1.9} and the property
$$
\chi_{\Delta^\prime } \cdot \chi_\Delta = \chi_{\Delta^\prime }
$$
implies
\begin{equation}
\label{ae1.31}
{\check{\Delta^\prime}} \subset {\check \Delta} \subset U.
\end{equation}
And the condition (i) implies that ${\check{\Delta^\prime}}$ is an open set.\\
In view of (iii) it follows that
$$
{\check \chi}_{\alpha_{g_i} ( \Delta^\prime )}\cdot
{\check \chi}_{\alpha_{g_j} ( \Delta^\prime )} = 0, \ \ \  i,j \in\overline{1,k}, \ i\neq j.
$$
Which means by the definition of $\alpha _g$ (see (\ref{ae1.15}) and
(\ref{ae1.16})), the interrelation between ${\hat T}_g$ and ${\check T}_g$ and
the definition of
$t_g$ for
${\check T}_g$ (see (\ref{e1.7}))
that
$$
t_{g_i} ({\check{\Delta^\prime}} ) \cap t_{g_j} ({\check{\Delta^\prime}} ) =
\emptyset \ \ \  i,j \in\overline{1,k}, \ i\neq j
$$
and consequently (in view of (\ref{ae1.31}))
the condition (\ref{e1.6'})
of the topologically free action is satisfied.\\

{\bf II.} Now let $\{ {\check T}_g \}_{g\in G}$ act topologically freely on $C(M)$.\\
Consider any finite set  $\{ g_1 , . . . , g_k \} \subset G$, any
element $\Delta \in \Sigma$ with $\mu (\Delta ) > 0$ and the set
$U= {\check \Delta} \subset M$ (as we have seen $U$ is a non empty open set).\\
Let $V\subset U$ be an open set mentioned in   (\ref{e1.6'}).
 Thus
\begin{equation}
\label{ae1.32}
t_{g_i} (V) \cap t_{g_j} (V) = \emptyset \ \ \  i,j \in\overline{1,k}, \ i\neq j.
\end{equation}
Take any $\Delta^\prime \in \Sigma$ satisfying the conditions
\begin{equation}
\label{ae1.33}
\mu ( \Delta^\prime ) >0
\end{equation}
 and
\begin{equation}
\label{ae1.34}
{\check {\Delta^\prime}} \subset V
\end{equation}
(such set does exist by (\ref{ae1.17})).\\
Observe that by the choice of $\Delta ^\prime$ we have
 $$
{\check \chi}_{\Delta}\cdot {\check \chi}_{\Delta^\prime} =
 \chi_U \cdot \chi_{\check{\Delta^\prime}} =
 {\check \chi}_{\Delta^\prime}
$$
which means that
$$
\chi_\Delta \cdot \chi_{\Delta^\prime} = \chi_{\Delta^\prime}
$$
that is
\begin{equation}
\label{ae1.35}
\Delta^\prime \subset \Delta \ ({\rm mod}\ \mu ).
\end{equation}
Now the properties (\ref{ae1.34}) and (\ref{ae1.32}) along with the argument
 used in {\bf I} shows that
$$
\mu (\alpha_{g_i} ( \Delta^\prime )\cap
\alpha_{g_j} ( \Delta^\prime ) )= 0, \ \ \  i,j \in\overline{1,k}, \ i\neq j
$$
which together with (\ref{ae1.33}) and (\ref{ae1.35}) proves the
metrical freedom
of the action of $\{ {\hat T}_g \}$.\\
Thus the 'coincidence' of the topological and metrical freedom is
established.
\mbox{}\qed\mbox{}\\

({\bf Remark.} We would like to note that certain interesting interrelations
between the algebras $L^\infty _\mu (\Omega )$ and $C(M)$ are described for example
in Rudin \cite{Rud}, 11.13, (f) ).\\

Now the analogue to
Theorem  \ref{1.7}
for the measurable case considered is
\begin{Tm}
\label{a1.11}
Let $A$ and $T_g , \ g\in G$ be those considered in \ref{a1.8}. If $G$
 acts metrically freely then $B(A, T_g )$ possesses the property (*).
\end{Tm}
{\bf Proof:} Just follow the idea of the proof of
Theorem 2.8 of \cite{Leb} (the statement in \cite{Leb} which corresponds to
Theorem \ref{1.7} of this paper)
%{1.7}%
using appropriate sets $\Delta , \ \Delta^\prime$ and the projection
$\chi_{\Delta^\prime }$ instead of the function $c(x)$ exploited in the proof of
Theorem 2.8.
\mbox{}\qed\mbox{}\\

\section{Example 1. Operators in $L^p (\Omega , E),   1<p<\infty $}
\setcounter{equation}{0}

\begin{I}
\label{a2.11}
\em
Let $(\Omega , \mu)$ be a space with a $\sigma -$additive
$\sigma -$finite measure $\mu$. Consider the space
$D = L^p_\mu (\Omega , E),\ \ 1<p<\infty $. Let
$A = L^\infty _\mu (\Omega , L(E)) \subset L(D)$ be the algebra of
{\em multiplication operators} defined by
\begin{equation}
\label{ae2.1}
(af)(x) = a(x)f(x), \ \ \ a\in A, \ f\in D
\end{equation}
and let $\alpha_g : \Omega \to \Omega , \ \ g\in G$ be a group of
measurable mappings preserving the equivalence class of $\mu$. By
$T_g$ we denote the isometry of $D$ defined by
\begin{equation}
\label{ae2.53}
(T_g f)(x) = \left[  \frac{d\alpha^{-1} _g (\mu )}{d\mu} \right]^{\frac{1}{p}}
f(\alpha^{-1} _g (x))
\end{equation}
where $\frac{d\alpha^{-1} _g (\mu )}{d\mu}$ is the Radon-Nikodim derivative
of $\alpha^{-1} _g (\mu )$ with respect to $\mu$.\\
One can easily verify that $T_g$ satisfies (\ref{ae1.1}) and the
mentioned mappings $\alpha_g$ coincide with those described in
\ref{a1.8}.\\

Let $B(A, T_g )\subset L(D)$ be the algebra generated by $A$ and
$\{ T_g \}_{g\in G}$.\\
\end{I}

The idea constituting the main plot of the investigation of this
example is similar to  that exploited in Example 2 in \cite{Leb}.
Therefore the consideration will consist of a series of steps.\\

The first one is
\begin{La}
\label{a2.12}
Let $B(A, T_g )$ be the algebra described in \ref{a2.11} and
$B({\bar A}, V_g )$ be the corresponding regular representation in the space
$H = l^p (G, L^p _\mu (\Omega , E))$. If $G$ acts metrically freely then the mapping
$$
B(A, T_g ) \to B({\bar A}, V_g )
$$
defined by
$$
a\to {\bar a}, \ \ \  a\in A
$$
$$
T_g \to V_g , \ \ \  g\in G
$$
is norm decreasing.
\end{La}
{\bf Proof:} The proof is in fact a certain inversion of the
argument used in the proofs of Lemmas 5.2 and 5.6 of \cite{Leb}.
%{2.7} and {2.9}%
\\

Consider any element ${\bar b} \in B({\bar A}, V_g )$ of the form
\begin{equation}
\label{ae2.54}
{\bar b} = \sum_{g\in F} {\bar a}_g V_g , \ \ \ \ \vert F \vert <\infty
\end{equation}
Fix any $\varepsilon >0$ and let $\eta = \{ \eta_g \} \in H$ be some vector
satisfying the conditions
\begin{equation}
\label{ae2.55}
\Vert \eta \Vert = 1
\end{equation}
and
\begin{equation}
\label{ae2.56}
\Vert {\bar b} \eta \Vert > \Vert {\bar b}  \Vert - \varepsilon .
\end{equation}
Without the loss of generality we can assume that
there exists a finite set $M\subset G$ such that
\begin{equation}
\label{ae2.57}
\eta_g =0 \ \ {\rm when}\ \ g\notin M
\end{equation}
Recalling the argument of the proof of
Lemma 5.6 \cite{Leb}
%{2.9}%
consider the vectors $\eta$ and ${\bar b}\eta$ belonging to the space
$L^p _\mu (\Omega , l^p (G, E))$.\\
Define the measure $\mu_\eta$ by the equation
\begin{equation}
\label{ae2.58}
\mu_\eta (\Delta ) = \int_\Delta \Vert \eta (x)\Vert^p \, d\mu
\end{equation}
for any measurable set $\Delta \subset \Omega$.\\
In view of (\ref{ae2.55}) we have
\begin{equation}
\label{ae2.59}
\mu_\eta (\Omega ) =1.
\end{equation}
now (\ref{ae2.59}) and (\ref{ae2.56})
imply
\begin{equation}
\label{ae2.60}
\sup_{\eta (x) \neq 0} \frac{\Vert ({\bar b}\eta )(x)  \Vert^p}{\Vert \eta (x)\Vert^p} >
(\Vert {\bar b} \Vert -\varepsilon)^p
\end{equation}
(note that by the construction ${\bar b}$ is a multiplication operator
in $L^p _\mu (\Omega , l^p (G, E))$ that is acting at $\eta$
according to formula (\ref{ae2.1})).\\
But (\ref{ae2.60}) means that there exists a $\mu_\eta$ measurable
(thus $\mu$ measurable) set $\Delta$ with $\mu_\eta (\Delta )>0$
(thus $\mu (\Delta )>0$) such that
\begin{equation}
\label{ae2.61}
 \frac{\Vert ({\bar b}\eta )(x)  \Vert}{\Vert \eta (x)\Vert} >
\Vert {\bar b} \Vert -\varepsilon \ \ \ {\rm for \ every}\ \ x\in \Delta
\end{equation}
Let
$$
F_1 = \cup_{g\in F} M\cdot g^{-1} .
$$
Since $G$ acts metrically freely it follows that there exists a measurable set
$\Delta^\prime \subset \Delta$ with $\mu (\Delta^\prime )>0$ and satisfying
the condition
\begin{equation}
\label{ae2.62}
\mu (\alpha^{-1} _{g_1} (\Delta^\prime ) \cap \alpha^{-1} _{g_2} (\Delta^\prime )) =0
\end{equation}
when $g_1 \neq g_2$ and $g_1 , g_2 \in [F_1 \cup M]$. \\
Consider the vector $\nu = \{ \nu_g \} \in H$ given by
$$
\nu_g = \chi_{\Delta^\prime}\cdot \eta_g
$$
As $\Delta^\prime \subset \Delta$ (\ref{ae2.61}) implies
\begin{equation}
\label{ae2.63}
\Vert {\bar b} \nu \Vert \le (\Vert {\bar b}  \Vert -\varepsilon ) \Vert \nu \Vert
\end{equation}
Now set the vector ${\bar \nu} \in L^p _\mu (\Omega , E)$ by
\begin{equation}
\label{ae2.64}
{\bar \nu}= \sum_{g\in M} \left(  T_{g^{-1}}\nu_g \right) =
\sum_{g\in M}{\bar \nu}_g .
\end{equation}
In view of (\ref{ae2.62}) we have
\begin{equation}
\label{ae2.65}
\Vert {\bar \nu} \Vert = \Vert \nu \Vert
\end{equation}
If $b\in B(A, T_g )$ is an operator of the form
\begin{equation}
\label{ae2.66}
b = \sum_{g\in F} a_g T_g
\end{equation}
then by the construction of ${\bar \nu}$ (see (\ref{ae2.64})) it follows that
\begin{equation}
\label{ae2.67}
(b {\bar \nu})(x) = 0 \ \ \  {\rm when} \ \ \
x\notin \cup_{g\in F_1}\alpha^{-1} _g  (\Delta^\prime )
\end{equation}
and for $g^\prime \in F_1$ we obtain
\begin{eqnarray*}
&&\left\Vert \chi_{\alpha^{-1} _{g^\prime}  (\Delta^\prime )}b{\bar \nu}\right\Vert =
\left\Vert \sum_{g\in F}a_g T_g {\bar \nu}_{g^\prime \cdot g}\right\Vert =
\left\Vert T_{g^\prime }
\left(\sum_{g\in F} a_g T_g T_{(g^\prime \cdot g)^{-1}} \nu_{g^\prime \cdot g}\right)
\right\Vert =\\
&&\left\Vert \sum_{g\in F} {\hat T}_{g^\prime}(a_g ) \nu_{g^\prime \cdot g} \right\Vert =
\left\Vert  \left(\left[ \sum_{g\in F} {\bar a}_g V_g \right] \nu\right)_{g^\prime}  \right\Vert=
\left\Vert ({\bar b}\nu)_{g^\prime} \right\Vert
\end{eqnarray*}
which along with (\ref{ae2.67}) and (\ref{ae2.62}) implies
$$
\Vert b{\bar \nu}\Vert = \Vert {\bar b}\nu \Vert
$$
and recalling (\ref{ae2.63}) and (\ref{ae2.65}) we conclude that
$$
\Vert b{\bar \nu}\Vert \ge (\Vert {\bar b}\Vert - \varepsilon) \Vert {\bar \nu}\Vert .
$$
Thus since $\varepsilon$ is arbitrary the proof is complete.
\mbox{}\qed\mbox{}\\

Now the analogue to  Lemma 5.5 \cite{Leb}
%{2.8}%
in the situation under consideration is
\begin{La}
\label{a2.13}
Let $B(A, T_g )$ be the algebra described in \ref{a2.11} and
$B({\bar A}, V_g)$ be the corresponding regular representation in the space
$H = l^p (G, L^p _\mu (\Omega , E))$. If $G$ is amenable then the mapping
$$
B({\bar A}, V_g) \to B(A, T_g )
$$
generated by the mappings
\begin{equation}
\label{ae2.68}
{\bar a} \to a, \ \ \ a\in A
\end{equation}
\begin{equation}
\label{ae2.69}
V_g \to T_g , \ \ \ g\in G
\end{equation}
is norm decreasing.
\end{La}
{\bf Proof:}
The same as for
Lemma 5.5 \cite{Leb}.
%{2.8}%
\mbox{}\qed\mbox{}\\

We can summarize the results obtained in
\begin{Tm}
\label{a2.14}
Let $B(A, T_g )$  and $B({\bar A}, V_g)$ be those considered in the lemma.\\
If $G$ is amenable and acts metrically freely then
$$
B(A, T_g )\cong    B({\bar A}, V_g)
$$
where the isomorphism is given by (\ref{ae2.68}) and (\ref{ae2.69}).\\
In particular the algebra $ B(A, T_g )$ possesses the properties (*) and
(**) and (\ref{e1.36}).
\end{Tm}
{\bf Proof:}
Apply Lemmas \ref{a2.12} and \ref{a2.13} and the results of
\ref{1.17}.
\mbox{}\qed\mbox{}\\

Recall that the consideration of
Example 2 \cite{Leb}
involved the interrelation between the operators
$b \in B(A, T_g ), \ \ {\bar b}\in B({\bar A}, V_g)$ and
their trajectorial representations
$b_x , \ x\in X$
(see Lemmas 5.2 and 5.6 \cite{Leb}),
%{2.7} and {2.9}%
where the latter representations are defined in the following way: \\
 for every point $x\in \Omega$ we define the representation
\begin{equation}
\label{ae2.24}
\pi_x : B(A, T_g ) \to L(l^p (G, E))
\end{equation}
$$
\pi_x (b) =b_x , \hspace{5mm}b\in B(A, T_g )
$$
by the equations
\begin{equation}
\label{ae2.25}
(\pi_x (a)\xi )_g = a(\alpha^{-1}_g (x))\xi_g ,
\end{equation}
\begin{equation}
\label{ae2.26}
(\pi_x (T_{g_0 })\xi )_g = \xi_{g g_0}
\end{equation}
where $\xi = (\xi_g  )_{g \in G} \in l^p (G, E)$ and $a \in A$.\\

We have not used the analogue of these interrelations in the
example considered yet.
And in order to obtain this analogue we have to impose certain
separability conditions on the measure $\mu$, group $G$ and space $E$
(this is because of the fact that the spaces of measurable functions
are rather 'large' and have 'more complicated' structure than the spaces of
continuous functions  (it is difficult to define what does it mean
the value of a measurable function at {\em a point})).\\

The next lemma presents the above mentioned analogue
in the measurable situation.
\begin{La}
\label{a2.15}
Let $(\Omega , \mu)$ be a separable measure space ($\mu$ be $\sigma -$finite),
$G$ be a countable group and $E$ be a separable Banach space.
Let $A$ and $\{ \alpha_g \}_{g\in G}$ be those described in \ref{a2.11}.
Then for any element
$$
{\bar b} = \sum_{g\in F}{\bar a}_g V_g \in B({\bar A}, V_g) , \ \ \ \vert F\vert <\infty
$$
we have
$$
\Vert {\bar b}\Vert = {\rm esssup}_\Omega  \Vert b_x \Vert
$$
where $b_x , \ x\in \Omega $ is defined by (\ref{ae2.24})-(\ref{ae2.26}).
\end{La}
{\bf Proof:}
The inequality
$$
\Vert {\bar b}\Vert \le {\rm esssup}_\Omega  \Vert b_x \Vert
$$
can be established by the argument absolutely analogous to that
used in the proof of
Lemma 5.6 \cite{Leb}.
%{2.9}%
So the point is to obtain the opposite inequality.\\

Since $(\Omega ,\mu )$ is separable we can assume without loss of
generality that it is strictly separable (see Halmos
\cite{Halmos}) that is there exists a sequence $\Delta_n$ with
$\mu (\Delta_n )<\infty \ \ n=1,2,...$ generating the $\sigma
-$algebra of all the measurable subsets of $\Omega $ (recall also
that the separability of $(\Omega ,\mu )$ is equivalent to the
separability of $L^p _\mu (\Omega ), \ \ 1\le p<\infty$).\\

We denote by ${\tilde D}$ the countable subset of
$$
l^p (G, L^p_\mu (\Omega , E))\cong L^p _{\mu \otimes \mu^\prime }(\Omega \times G, E)
\cong L^p _\mu (\Omega , l^p (G,E))
$$
($\mu^\prime $ is the discrete measure on $G$) consisting of finite
sums of vectors of the form
$$
\chi (\Delta_n \times g_m )e_k \ \ \ n,m,k = 1,2,...
$$
where $g_m \in G, \ \ \{ e_k \}$ is a countable dense subset of
$E, \ \ \chi (\Delta_n \times g_m )$ is the characteristic function of the set
$\Delta_n \times g_m \subset \Omega \times G$.\\
The definition of ${\tilde D}$ implies that the set of vectors
$$
\{  \xi (x) \}_{\xi \in {\tilde D} }\subset l^p (G, E)
$$
is dense in $ l^p (G, E)$ for every $x\in \Omega$.\\
Now for any $f\in L^\infty _\mu (\Omega )$ and for any
$\eta = \{ \eta (x) \}_{x\in \Omega }\in L^p _\mu (\Omega ,  l^p (G, E))$ we have
$$
f\cdot \eta \in  L^p _\mu (\Omega ,  l^p (G, E))
$$
and
\begin{eqnarray*}
&&\int_\Omega \vert f(x) \vert^p  \Vert ({\bar b}\eta )(x) \Vert^p \, d\mu =
 \Vert {\bar b}(f\eta )(x) \Vert^p \le
\Vert {\bar b}\Vert^p   \Vert f\eta \Vert^p =\\
&& \Vert {\bar b}\Vert^p \int_\Omega \vert f(x) \vert^p  \Vert \eta (x) \Vert^p \, d\mu
\end{eqnarray*}
which as $f$ is arbitrary implies
$$
\Vert b_x \eta (x) \Vert = \Vert ({\bar b}\eta)(x)\Vert \le \Vert {\bar b}\Vert
\Vert \eta (x) \Vert
$$
for almost every $x$.\\
Thus applying this inequality to vectors
$\xi \in {\tilde D}$ we have
$$
\Vert b_x \xi (x) \Vert \le \Vert {\bar b}\Vert \Vert \xi (x) \Vert \ \ \
{\rm almost \ \ everywhere}
$$
and since ${\tilde D}$ is a countable set and
$\{ \xi (x) \}_{\xi \in {\tilde D} }$ is dense in $l^p (G,E)$ it follows that
$$
{\rm esssup}_\Omega \Vert b_x \Vert \le \Vert {\bar b}\Vert
$$
and the proof is complete.
\mbox{}\qed\mbox{}\\

  \section{Example 2. Operators in $L^\infty _\mu (\Omega , E)$ }
\setcounter{equation}{0}
\begin{I}
\label{a2.20}
\em
Now we would like to present a measurable analogue to
Example  3 \cite{Leb}.

Let $D= L^\infty _\mu (\Omega , E)$ (where $\Omega $ is
 a space with a $\sigma -$additive $\sigma -$finite  measure $\mu$)
and $A = L^\infty _\mu (\Omega , L(E))$ be the algebra of multiplication operators
(see (\ref{ae2.1})) and
$\alpha_g : \Omega \to \Omega , \ g\in G$ be a group of measurable
mappings preserving the  equivalence class of $\mu $.
By $T_g$ we denote an isometry of $D$ given by
\begin{equation}
\label{ae2.71}
(T_g f)(x) = f (\alpha^{-1}_g (x)) .
\end{equation}
Let $B(A, T_g )$   be the algebra generated by $A$ and $\{  T_g \}_{g\in G}$.\\

For any fixed finite set $F\subset G$  we denote by $B_F (D)$ and $S_F (D)$
respectively the sets
\begin{equation}
\label{ae2.3}
B_F (D) = \{  \{  f_g  \}_{g\in F} : f_g \in D, \hspace{2mm}\Vert f_g \Vert \le 1,  \hspace{2mm}g\in F  \} ,
\end{equation}
\begin{equation}
\label{ae2.4}
S_F (D) = \{  \{  f_g  \}_{g\in F} : f_g \in D, \hspace{2mm}\Vert f_g \Vert = 1,  \hspace{2mm}g\in F  \} .
\end{equation}
We begin with the statement which asserts a bit less than
Theorem 6.2 \cite{Leb}
% {2.17}%
(from the 'measurable' point of view) but anyway it is equivalent to it
(see Remarks \ref{a2.22} \ (2), (3)).
\end{I}
\begin{Tm}
\label{a2.21}
Let $B(A, T_g )$  be the  algebra  introduced above. If $G$ acts
metrically freely then for any finite
$F$ we have
\begin{equation}
\label{ae2.72}
\Vert \sum_{g\in F} a_g T_g \Vert =
\sup_{ \{ f_g \}_{g\in F} \in S_F (D)} \Vert  \sum_{g\in F} a_g f_g \Vert =
\sup_{\{ f_g \}_{g\in F} \in B_F (D)} \Vert  \sum_{g\in F} a_g f_g \Vert
\end{equation}
where $B_F (D)$ and $S_F (D)$ are defined by (\ref{ae2.3}) and (\ref{ae2.4}).
\end{Tm}
{\bf Proof:}
The proof is in fact a measurable variant of the proof of
Theorem 4.2 \cite{Leb}.\\
%{2.2}%
We begin with the establishing of the equality
\begin{equation}
\label{ae2.73}
\Vert \sum_{g\in F} a_g T_g \Vert =
\sup_{ \{ f_g \}_{g\in F} \in B_F (D)} \Vert  \sum_{g\in F} a_g f_g \Vert
\end{equation}
First  of all for any $f \in L^\infty _\mu (\Omega , E) = D$ with
$\Vert f \Vert =1$ we have
$$
\Vert \sum_ F a_g (T_g f) \Vert =  \Vert \sum_ F a_g f_g \Vert \le
\sup_{ \{ f_g \}_{g\in  F} \in B_F (D)} \Vert  \sum_ F a_g f_g \Vert
$$
which means that
$$
\Vert \sum_ F a_g T_g \Vert \le
\sup_{ \{ f_g \}_{g\in F} \in B_F (D)} \Vert  \sum_F a_g f_g \Vert
$$
To establish the opposite inequality fix any collection
${ \{ f_g \} \in B_F (D)} $ and let
\begin{equation}
\label{ae2.74}
\Vert \sum_ F a_g f_g \Vert  = \lambda
\end{equation}
For any $\varepsilon >0$ (\ref{ae2.74}) implies the existence of a
measurable set $\Delta $ with $\mu (\Delta )> 0$ and such that
\begin{equation}
\label{ae2.75}
\Vert \sum_ F a_g (x) f_g (x) \Vert  = \lambda -\varepsilon \ \
{\rm for \ \ every}\ \ x\in \Delta
\end{equation}
Since $G$ acts metrically freely it follows that there exists a measurable set
$\Delta^\prime \subset \Delta , \ \ \mu (\Delta ^\prime >0)$ satisfying the condition
\begin{equation}
\label{ae2.76}
\mu (\Delta ^\prime \cap \alpha^{-1}_g (\Delta ^\prime )) = 0, \ \ \ g\in F,\ g\neq e
\end{equation}
Consider the function
\begin{equation}
\label{ae2.77}
f = \sum_F \chi_{\alpha^{-1}_g (\Delta ^\prime )}f_g
\end{equation}
As $\{  f_g \} \in B_F (D)$ we have
$$
\Vert f \Vert \le 1
$$
and on the other hand by the explicit form of $f$ and (\ref{ae2.75})
we conclude that
$$
{\rm for \ \ every}\ \ x\in \Delta^\prime \ \ \
\Vert (\sum_ F a_g T_g f) (x) \Vert  = \lambda -\varepsilon
$$
which by the definition of $\lambda$ (see (\ref{ae2.74})) and in
view of the arbitrariness of $\varepsilon$ means that
$$
\Vert \sum_{g\in F} a_g T_g \Vert \ge
\sup_{ \{ f_g \} \in B_F (D)} \Vert  \sum_{g\in F} a_g f_g \Vert
$$
thus finishing the verification of (\ref{ae2.73}).\\

The equivalence of (\ref{ae2.73}) and (\ref{ae2.72}) is established by
the same argument as the one used for the corresponding purposes
in the proof of Theorem 4.2 \cite{Leb}.\\
% {2.2}%
The proof is complete.
\mbox{}\qed\mbox{}\\
\begin{Rk}
\label{a2.22}
\em \ \\

(1) If $E = {\mathbb C}$ (that is $D= L^\infty _\mu (\Omega )$ and
 $A = L^\infty _\mu (\Omega )$) then (\ref{ae2.72}) implies \\

{\em
if $G$ acts metrically freely then
$$
\Vert \sum_F a_g T_g \Vert = \ \ {\rm esssup}_\Omega \sum_F \vert a_g (x) \vert
$$
}

Indeed on the one hand
$$
\sup_{ \{ f_g \} \in B_F (D)} \Vert  \sum_{g\in F} a_g f_g \Vert \le
 {\rm esssup}_\Omega \sum_F \vert a_g (x) \vert
$$
and to obtain the opposite inequality just set
$$
f_g (x) = \left\{
\begin{array}{ccl}
[{\rm arg} \, a_g (x)]^{-1}&,& \  {\rm if}\ \ a_g (x)\neq 0\\
1&,& \  {\rm if}\ \ a_g (x) = 0
\end{array}
\right.
$$
(since $a_g \in L^\infty _\mu (\Omega ) \ \ \ f_g \in L^\infty _\mu (\Omega )$
as well).\\

(2) The equality (\ref{ae2.72}) also shows that\\

{\em
if $G$ acts metrically freely then
$$
\Vert  \sum_F a_g T_g  \Vert = \Vert {\tilde b}_F \Vert
$$
 where
$$
 {\tilde b}_F : D_1 \times D_2 \times ... \times D_{\vert F \vert} \to D, \hspace{5mm} D_i = D
$$
 is given by
\begin{equation}
\label{ae2.22}
{\tilde b}_F (\xi_1 , ... , \xi_{\vert F \vert}) = a_{g_1}\xi_1 + ... +
 a{g_{\vert F \vert} }\xi_{\vert F \vert}
\end{equation}
($\{ g_1 , ... , g_{\vert F \vert}  \} = F$).
} \\

(3) The preceding remark leads in turn to the next observation.\\

Since $a_g \in L^\infty _\mu (\Omega , L(E)), \  g\in F$ it follows (by the structure of
${\tilde b}_F$) that
$$
{\tilde b}_F \in  L^\infty _\mu (\Omega , L({\tilde E},E))
$$
where ${\tilde E} = E_1 \times E_2 \times . . . \times E_{\vert F \vert }, \ \ E_i = E$.\\
But this means that
$$
\Vert {\tilde b}_F (\cdot )  \Vert \in L^\infty _\mu (\Omega )
$$
and
$$
\Vert {\tilde b}_F   \Vert = \ \ {\rm esssup}_\Omega \Vert {\tilde b}_F (x )  \Vert
$$
which along with the preceding remark (2) implies\\

{\em if $G$ acts metrically freely then}
\begin{eqnarray}
\label{ae2.78}
&& \Vert \sum_F a_g T_g \Vert = \ \ {\rm esssup}_\Omega \
\sup_{ \{ f_g \} \in B_F (E)} \Vert  \sum_{g\in F} a_g (x) f_g \Vert = \nonumber \\
&& {\rm esssup}_\Omega \  \sup_{ \{ f_g \} \in S_F (E)} \Vert  \sum_{g\in F} a_g (x) f_g \Vert
\end{eqnarray}
thus strengthening the statement of Lemma \ref{a2.21} up to the
 measurable variant of Theorem 6.2 \cite{Leb}.
% {2.17}%
\end{Rk}

Now we finish the consideration of the example with the measurable
analogue to
Lemma 6.4 \cite{Leb}.
% {2.19}%
\begin{La}
\label {a2.23}
Let  $B(A, T_g )$     be the algebra described in \ref{a2.20} and
 $B({\bar A}, V_g)$ be the corresponding regular representation in the space
$$
H = l_0 (G, L^\infty _\mu (\Omega , E)) \ \ ({\rm or}\ \  l^\infty  (G, L^\infty _\mu (\Omega , E)) ).
$$

If $G$ acts metrically freely then
$$
B(A, T_g )\cong B({\bar A}, V_g)
$$
where the isomorphism is generated by the mappings
$$
a\to {\bar a}, \ \ \ a\in A
$$
$$
T_g \to V_g , \ \ \ g\in G
$$
and in particular $B(A, T_g )$ possesses the properties (*) and (**) and
(\ref{e1.36}).
\end {La}
{\bf Proof:}
Simple calculation shows that the norm of the element
$$
{\bar b} = \sum_{g\in F}{\bar a}_g V_g \in B({\bar A}, V_g)
$$
is equal to
$$
\sup_{ \{ f_g \} \in B_F (D)} \Vert  \sum_{g\in F} a_g  f_g \Vert
$$
which establishes the isomorphism
$$
B(A, T_g )\cong B({\bar A}, V_g).
$$
And consequently the properties (*), (**) and
(\ref{e1.36})
for  $B(A, T_g )$ follow from the results of
 \ref{1.17}.
\mbox{}\qed\mbox{}\\

\section{Example 3. Operators in $L^1_\mu (\Omega , E)$}
\setcounter{equation}{0}

Our final example here deals with the space $L^1$ and in fact
 is a measurable variant of
Example 4 \cite{Leb}. The argument in the situation examined
hereafter is similar to that exploited  when considering Example 4
\cite{Leb} (that is we will reduce this case to the $L^\infty $
situation already studied) but  requires the usage of the results
of the preceding example instead of the results
of Example 3 \cite{Leb}.\\
\begin{I}
\label{a2.29}
\em
Let $(\Omega ,\mu )$ be the space considered in \ref{a2.11},
$D = L^1 _\mu (\Omega , E )$ and $A = L^\infty  _\mu (\Omega , L(E))$
be the algebra of operators defined by (\ref{ae2.1}) and $T_g$ be
defined by (\ref{ae2.53}) (with $p=1$).\\
Let   $B(A, T_g ) \subset L(D)$ be the algebra generated by
$A$ and $\{ T_g \}_{g\in G}$.
\end{I}
\begin{I}
\label{a2.30} \em Bearing in mind the reasoning of Example 4
\cite{Leb} (that is introducing and using the corresponding
formally adjoint operators) one can obtain the next statement
(which is the natural substitute of the result presented in 7.3
\cite{Leb}
% {2.26}%
 for the operators considered).\\

{\em
Let $B(A, T_g )$ be the algebra introduced in \ref{a2.29}.
If $G$ acts metrically freely then}
\begin{eqnarray}
\label{ae2.88}
&& \Vert \sum_{g\in F} a_g T_g \Vert = \ \ {\rm esssup}_\Omega \
\sup_{ \{ f_g \} \in B_F (E^\ast )} \Vert
 \sum_{g\in F} [a_g (\alpha_g (x))]^\ast f_g \Vert = \nonumber \\
&& {\rm esssup}_\Omega \  \sup_{ \{ f_g \} \in S_F (E^\ast )} \Vert
 \sum_{g\in F} [a_g (\alpha (x))]^\ast  f_g \Vert
\end{eqnarray}
{\bf Proof:}
Just follow the argument of
7.3 \cite{Leb}
%{2.26}%
using Remark \ref{a2.22} (3) instead of
Theorem 6.2 \cite{Leb}.
%{2.17}%
\mbox{}\qed\mbox{}\\
\begin{Rk}
\label{a2.31}
\em
(see Remark \ref{a2.22} (1)).\\
If $E = {\mathbb C}$ and $G$ acts metrically freely then
$$
\Vert \sum_F a_g T_g \Vert = \ \ {\rm esssup}_\Omega \
 \sum_F \vert a_g (\alpha_g (x))\vert .
$$
\end{Rk}
And finally the analogue to  Lemma \ref{a2.23}
here is
\begin{La}
\label{a2.32}
Let $B(A, T_g )$  be the algebra described in \ref{a2.29} and
$B({\bar A}, V_g)$  be the corresponding regular representation in the space
$l^1 (G, L^1 _\mu (\Omega , E))$.
If $G$ acts metrically freely then
$$
B(A, T_g )\cong B({\bar A}, V_g)
$$
where the isomorphism is generated by the mappings
$$
a\to {\bar a}, \ \ \ a\in A
$$
$$
T_g \to V_g , \ \ \ g\in G
$$
and in particular $B(A, T_g )$ possesses the properties (*) and (**) and
(\ref{e1.36}).
\end{La}
\end{I}

\section{Isomorphism theorems. Interpolation}
\setcounter{equation}{0}
Now we would like to observe certain interrelations
between the examples considered and the
Isomorphism Theorem (\cite{AnLeb}, Corollary 12.17).
 \begin{I}
\label{a2.33}
\em
(1) Observe that in Examples 2, 3 we did not use any information about
the group $G$ thus the group in these examples is
{\em not necessarily amenable}.\\

(2) The essentially different picture is drawn in Example 1.\\
Here\\
(i) if $G$ acts metrically freely then \\
$B({\bar A}, V_g )$ is a representation of $B(A, T_g )$ for {\em any} $G$
({\em not necessarily amenable})\\
 (Lemma \ref{a2.12}).\\

While\\
(ii) if $G$ is amenable then\\
 $B(A, T_g )$ is a representation of $B({\bar A}, V_g )$ for an {\em arbitrary} action of $G$
({\em not necessarily metrically free}).\\
(Lemma \ref{a2.13}).\\

Thus in these examples the metrical freedom of the action of $G$
and the amenability
of $G$ are lying in a sense opposite each other.\\
If $G$ acts metrically freely then  $B(A, T_g )$  is 'larger' than
$B({\bar A}, V_g )$
(see (i)).\\
And\\
if $G$ is amenable then  $B({\bar A}, V_g )$ is 'larger' than  $B(A, T_g )$ (see (ii)).\\

Both these algebras 'coincide' if $G$ acts metrically freely and
is amenable
(Theorem \ref{a2.14}).\\

(3) Consideration of Example 1 leads to certain {\em Isomorphism Theorems} which
(just as it was done in  \cite{AnLeb}, Corollary 12.17) establish the isomorphism between {\em essentially
spatially different} operator algebras (thus wiping off the spaces where these operators  act).\\
For example.\\
Let $(\Omega , \mu_i ), \ \ i=1,2$ be two spaces with $\sigma -$additive
$\sigma -$finite separable measures $\mu_1$ and $\mu_2$ absolutely continuous
with respect to each other (thus
$L^\infty _{\mu_1} (\Omega, L(E)) \cong L^\infty _{\mu_2} (\Omega, L(E))$)
and let $\{ \alpha_g \}_{g\in G}$ be a group of measurable mappings of $\Omega$
preserving the equivalence classes of $\mu_1$ and $\mu_2$.
Consider the spaces $D_i = L^p _{\mu_i} (\Omega, E), \ \i=1,2$.
 Let $A_i = L^\infty _{\mu_1} (\Omega, L(E)) \subset L(D_i )$ be the
algebras of multiplication operators defined by (\ref{ae2.1}) and $T^i _g , \ \ i=1,2$
be the isometries of $D_i$ defined by (\ref{ae2.53}) (with $\mu = \mu_i$) and
$B(A, T^i _g )$ be the algebras generated by $A_i$ and $\{ T^i _g \}_{g\in G}$.\\
The Isomorphism Theorem related to Example 1 is stated as follows:\\

{\em If $E$ is a separable Banach space and $G$ is a countable
amenable group acting metrically freely then $B(A_1 , T^1 _g )$
and  $B(A_2 , T^2 _g )$ are isomorphic (as Banach algebras) and
the isomorphism is established by the natural isomorphism
$$
A_1 \cong A_2
$$
and the mapping
$$
T^{1}_g \to T^{2}_g .
$$
 }
{\bf Proof:}
Follows from Lemma \ref{a2.15}.
\mbox{}\qed\mbox{}\\
\end{I}

We have not obtained the explicit formula for the norm of operator
$\sum a_g T_g$ in Example 1. But in fact the formulae proved while
considering the other examples
and the Riesz-Thorin interpolation theorem makes it possible to write out the useful
estimates for the norm of $\sum a_g T_g$ in all the cases (in a way analogous
to that presented  in Section 8 of \cite{Leb}.\\
\begin{I}
\label{a2.39}
\em
Let $b \in B(A, T_g ) \subset L(L^{p}_\mu   (\Omega ,E)), \ \  1\le p\le \infty $
be an operator of the form
\begin{equation}
\label{ae2.91}
b = \sum_{g\in F}a_g T_g , \hspace{5mm} \vert F \vert < \infty
\end{equation}
where for $1\le p < \infty $ $B(A, T_g )$ is that described in \ref{a2.11}
and for $p = \infty$ in \ref{a2.20}.\\
We shall denote by
$\Vert b\Vert_p$ the norm of the operator $b$ as acting in the space
$L^{p}_\mu   (\Omega ,E))$. \\
For these objects the measurable substitution of Lemma 8.2 \cite{Leb}
%2.37%
is
\begin{La}
\label{a2.40}
If $G$ is amenable and acts metrically freely on $\Omega$ then for any
$p\in (1,\infty )$ we have
\begin{equation}
\label{ae2.99}
\Vert b \Vert_p \le \Vert b \Vert^{1\over p}_1 \cdot \Vert b \Vert^{1 - {1\over p}}_\infty
\end{equation}
where $\Vert b \Vert_1$ is given by (\ref{ae2.88}) and  $\Vert b \Vert_\infty$ by (\ref{ae2.78}).
\end{La}
\end{I}


\begin{thebibliography}{                                  }

\bibitem{Leb}
A. Lebedev,
{\em  On the structure of Banach algebras associated with automorphisms.}
(To appear).

\bibitem{Ped}
G.K. Pederesen, {\em $C^\ast -$algebras and their automorphisms groups.}
Academic Press, 1979.

\bibitem{Land}
M.B. Landstad, Duality theory for covariant systems, {\em Trans. Amer. Math. Soc.,}
1979, V. 248, p. 223-267.



\bibitem{AnLeb}
A. Antonevich, A. Lebedev, {\em Functional differential equations: I. $C^\ast -$theory.}
Longman Scientific $\&$ Technical, 1994.


\bibitem{DunSchw}
N. Dunford, J.T Schwarts, {\em Linear operators. Part 1: General theory}
Interscience, 1958.

\bibitem{Rud}
W. Rudin, {\em Functional analysis} McGraw-Hill, 1973.

\bibitem{Halmos}
P.R. Halmos, The decomposition of measures, {\em Duke Math. J.,}
1941, N 8, p. 386-392.





\end{thebibliography}
\end{document}